\documentclass[12pt]{amsart}

\usepackage{fullpage}
\usepackage{bbm}

\newcommand\g{\mathfrak g}
\renewcommand\l{\mathfrak l}
\newcommand\m{\mathfrak m}
\newcommand{\p}{\mathfrak{p}}%
\renewcommand{\t}{\mathfrak{t}}%
\renewcommand{\sl}{\mathfrak{sl}}%

\newcommand\cI{\mathcal I}
\newcommand\cVA{\mathcal{VA}}

\newcommand\va{\mathbf a}
\newcommand\vTheta{\mathbf \Theta}

\newcommand\bC{\mathbb C}
\newcommand\bF{\mathbb F}
\newcommand\bQ{\mathbb Q}
\newcommand\bZ{\mathbb Z}
\newcommand\bk{\mathbbm k}

\newcommand\sub{\subseteq}
\newcommand\op{\mathrm{op}}
\newcommand\into{\hookrightarrow}

\DeclareMathOperator{\ad}{ad}
\DeclareMathOperator{\End}{End} %
\DeclareMathOperator{\gr}{gr} %

\numberwithin{equation}{section}

\newtheorem{thm}[equation]{Theorem}
\newtheorem{lem}[equation]{Lemma}
\newtheorem{prop}[equation]{Proposition}

\theoremstyle{remark}
\newtheorem{rem}[equation]{Remark}

\title%[Representations of finite $W$-algebras]
{On $1$-dimensional
representations of finite $W$-algebras associated to simple
Lie algebras of exceptional type}

\author{Simon M.~Goodwin, Gerhard R\"ohrle and Glenn Ubly}
\address{School of Mathematics, University of Birmingham, Birmingham, B15
2TT, UK}  \email{goodwin@maths.bham.ac.uk}
\address{Fakult\"at f\"ur Mathematik, Ruhr-Universit\"at Bochum,
D-44780 Bochum, Germany} \email{gerhard.roehrle@rub.de}
\address{School of Mathematics, University of Southampton,
Southampton, SO17 1BJ, UK} \email{gu@soton.ac.uk}
 


\begin{document}

\begin{abstract}
We consider the finite $W$-algebra $U(\g,e)$ associated to a nilpotent element
$e \in \g$ in a simple complex Lie algebra $\g$ of exceptional type.
Using presentations obtained through an algorithm based on the PBW-theorem,
we verify a conjecture of Premet, that
$U(\g,e)$ always has a $1$-dimensional representation, when
$\g$ is of type $G_2$, $F_4$, $E_6$ or $E_7$.  Thanks to a
theorem of Premet, this allows one to deduce the
existence of minimal dimension representations of reduced enveloping
algebras of modular Lie algebras of the above types.  In addition,
a theorem of Losev allows us to
deduce that there exists a completely prime
primitive ideal in $U(\g)$ whose associated variety is the coadjoint orbit
corresponding to $e$.
\end{abstract}

\maketitle

\section{Introduction}
Finite $W$-algebras were introduced into the mathematical literature
by Premet in \cite{Pr2}.  We recall that a finite $W$-algebra
$U(\g,e)$ is a certain finitely generated algebra associated to a
reductive Lie algebra $\g$ and a nilpotent element $e \in \g$; it is
a quantization of the Slodowy slice through the nilpotent orbit of
$e$. Recently, there has been much research activity in the
representation theory of $U(\g,e)$, see for example \cite{BK2, BGK,
Gi, Lo1, Lo2, Lo3, Pr3, Pr4, Pr5}. However, the fundamental question
of the existence of a 1-dimensional representation of $U(\g,e)$
remains open in general.  Indeed it was only recently proved by
Premet in \cite[Cor.\ 1.1]{Pr4} that $U(\g,e)$ always has
finite-dimensional representations; alternative proofs have been
given independently by Losev \cite[Thm.\ 1.2.2(viii)]{Lo1} and
Ginzburg \cite[Thm.\ 4.5.2]{Gi}.  It was conjectured by Premet in
\cite[Conj.\ 3.1]{Pr3} that there is always a $1$-dimensional
$U(\g,e)$-module; this conjecture has now been verified for $\g$ of
classical type by Losev in \cite[Thm.\ 1.2.3]{Lo1}.  Further, the
conjecture was reduced to the case where $e$ is rigid by Premet in
\cite[Thm.\ 1.1]{Pr5}; we recall that $e$ is {\em rigid} if the
nilpotent orbit of $e$ cannot be obtained through Lusztig--Spaltenstein
induction \cite{LS} from a nilpotent orbit in a Levi subalgebra of
$\g$.

In this article we consider Premet's conjecture in case $\g$ is of
exceptional type: our main result is the following theorem.

\begin{thm} \label{T:main}
Let $\g$ be a simple Lie algebra over $\bC$ of type $G_2$, $F_4$, $E_6$ or
$E_7$, and let $e \in \g$ be nilpotent.  Then there exists a $1$-dimensional
representation of the finite $W$-algebra $U(\g,e)$.
\end{thm}

Our approach is computational and based on the PBW-theorem for
finite $W$-algebras, see \cite[Thm.\ 4.6]{Pr2} along with the above
mentioned theorem of Premet reducing to the case where $e$ is rigid.
We have developed an algorithm that determines an explicit
presentation of a given finite $W$-algebra, which we have used for
the cases $\g$ of type $G_2$, $F_4$, $E_6$ or $E_7$ and $e$ rigid. From
these presentations, it is straightforward to determine all
$1$-dimensional representations of $U(\g,e)$; in fact, the number of
1-dimensional representations is either one or two, see Table
\ref{Tab:rigidres}.  We remark that it is possible to determine the
$1$-dimensional representations without calculating a full
presentation, which we have done in a number of cases when $\g$ is
of type $E_7$, see Proposition \ref{P:1-dim} and Section
\ref{S:results} for more details.

Our results taken with the work of Losev and Premet
mean that Premet's  conjecture is reduced to the cases where $\g$
is of type $E_8$ and $e$ lies in a rigid orbit of $\g$.  Our methods
allow us to prove the
existence of $1$-dimensional $U(\g,e)$-modules for some rigid $e$ in case
$\g$ is of type $E_8$; at present it is computationally unfeasible
to deal with those $e$ with large height, see Remark \ref{R:E8}
for more details.\footnote{
After this paper was written,
I.\ Losev announced in \cite{Lo4} an alternative approach to finding
the $1$-dimensional
representations of $U(\g,e)$ based on the highest weight theory
from \cite{BGK}.}

%A notable observation is that the presentations that we have calculated,
%have coefficients only in $\bZ[d^{-1}]$, where $d$ is the product of
%the bad primes for $\g$.  This tells us that $U(\g,e)$ is defined
%over $\bZ[d^{-1}]$.  We also note here that this is the case for
%presentations of $U(\g,e)$ for non-rigid $e$ that we have calculated.

\smallskip

We now discuss how results of Premet provide an
application of Theorem \ref{T:main} to the representation
theory of modular Lie algebras, which we state in Theorem \ref{T:mindim}.
It is necessary to introduce some notation. % for its statement.

Let $G$ be a simple simply connected algebraic group over $\bC$. Let
$\g$ be the Lie algebra of $G$ and let $\g_\bZ$ be a Chevalley
$\bZ$-form of $\g$.  Let $p$ be a prime that is assumed to be good
for $G$ and let $\bk$ be the algebraic closure of $\bF_p$. Let
$\g_\bk = \g_\bZ \otimes \bk$ and let $\xi \in \g_\bk^*$. The
reduced enveloping algebra corresponding to $\xi$ is denoted
$U_\xi(\g_\bk)$.  We write $d_\xi$ for half the dimension of the
coadjoint orbit of $\xi$.

Given a nilpotent orbit in $\g$ it is possible to choose a
representative $e \in \g_\bZ$ such that we can find an
$\sl_2$-triple $(e,h,f)$ contained in $\g_\bZ$, see \S
\ref{ss:triple}. Let $\kappa$ be the Killing form on $\g$ and assume
that $p$ is such that it does not divide $\kappa(e,f)$.  Then we can normalize
$\kappa$ to obtain the bilinear form $(\cdot\,, \cdot)$ with $(e,f)
= 1$, which can be viewed as a bilinear form on both $\g$ and
$\g_\bk$. We may view $e$ both as an element of $\g$ and $\g_\bk$,
and therefore view $\chi$ defined by $\chi(x) = (x,e)$ both as an
element of $\g^*$ and $\g_\bk^*$.

The Kac--Weisfeiler conjecture, proved by Premet in \cite{Pr1}, says
that the dimension of a $U_\xi(\g_\bk)$-module is divisible by
$p^{d_\xi}$. The existence of a $U_\xi(\g_\bk)$-module with
dimension equal to $p^{d_\xi}$ is an open problem, see \cite[\S
4.4]{Pr1}. Thanks to the Kac--Weisfeiler theorem, or its
generalization \cite[Thm.\ 3.2]{FP} due to Friedlander and Parshall,
there is a reduction to the case of nilpotent $\xi$.

In \cite[Thm.\ 1.4]{Pr5}, Premet proves that the existence
of a $1$-dimensional $U(\g,e)$-module implies existence of a
representation of $U_\chi(\g_\bk)$ of dimension
$p^{d_\chi}$ for $p$ sufficiently large.  As a consequence of Premet's and
Losev's results, the existence of a minimal dimension representation of
$U_\xi(\g_\bk)$ is now proved for $\g$ of classical type and $p$ sufficiently
large, see \cite[Cor.\ 1.1]{Pr5}.  We can deduce the following
from Theorem \ref{T:main}.

\begin{thm} \label{T:mindim}
Let $\g_\bk$ be a simple Lie algebra over $\bk$ of type $G_2$,
$F_4$, $E_6$ or $E_7$. Assume that $p \gg 0$ and let $\xi \in
\g_\bk^*$.  Then the reduced enveloping algebra $U_\xi(\g_\bk)$ has
a simple module of dimension $p^{d_{\xi}}$, where $d_\xi$ is half
the dimension of the coadjoint orbit of $\xi$.
\end{thm}

The restriction in \cite[Thm.\ 1.4]{Pr5} leading to the condition
that $p$ is sufficiently large stems from a number of rationality
assumptions in the construction of $U(\g,e)$ made in \cite[\S
2]{Pr5}; in particular, that $U(\g,e)$ is defined over
$\bZ[d^{-1}]$, but only for $d$ a sufficiently large integer.  As
discussed in \cite[Rem.\ 2.2]{Pr5}, this bound can be lowered by
knowledge of explicit presentations of $U(\g,e)$.  In the cases
where $\g$ is of type $G_2$, $F_4$ and $E_6$, we are able to give an
explicit bound, which just requires $p$ to be good and not to divide
all $\kappa(e_0,f_0)$, where $(e_0,h_0,f_0)$ is an $\sl_2$-triple in
a Levi subalgebra of $\g$ in which $e_0$ is rigid; see Remark
\ref{R:char} for more details.

\smallskip

Next we discuss how Theorem \ref{T:main} implies existence of completely prime
primitive ideals, through Skryabin's equivalence and a result of
Losev.  This is stated in Theorem \ref{T:comprime} below.

We recall that Skryabin's equivalence gives an equivalence of categories
between the category of $U(\g,e)$-modules and a certain category of
generalized Whittaker modules, see \cite{Sk} or \cite[Thm.\ 6.1]{GG}.
For a $U(\g,e)$-module $M$, we write $\widehat M$ for the Whittaker
$U(\g)$-module obtained through this equivalence, and $\cI_M$
for the annihilator of $\widehat M$ in $U(\g)$.  Premet proved in
\cite[Thm.\ 3.1]{Pr2}
that the associated variety $\cVA(\cI_M)$ of $\cI_M$
contains the coadjoint orbit $G \cdot \chi$ and that
$\cVA(\cI_M) = G \cdot \chi$ if and only if $M$ is finite-dimensional.
Further, it was conjectured in \cite[Conj.\ 3.1]{Pr3}
that if $M$ is $1$-dimensional, then $\cI_M$ is completely prime;
this conjecture was verified by Losev in \cite[Prop.\ 3.4.6]{Lo1}.
Therefore, we can deduce the following from
Theorem \ref{T:main}.

\begin{thm} \label{T:comprime}
Let $\g$ be a simple Lie algebra over $\bC$
of type $G_2$, $F_4$, $E_6$ or $E_7$.
Let $e \in \g$ be nilpotent and let $\chi = (e,\cdot) \in \g^*$.
Then there is a completely prime
primitive ideal of $U(\g)$ whose associated variety is $G \cdot \chi$.
\end{thm}

We give a brief outline of the structure of the paper.
First we recall the definition of the finite $W$-algebra $U(\g,e)$ and
state a version of the PBW-theorem for $U(\g,e)$
in Section \ref{S:prelim}.  Next in Section \ref{S:small}, we use
the PBW-theorem to give presentations of $U(\g,e)$ and show,
in Theorem \ref{T:small}, that some
of the relations are superfluous.
Our algorithm is explained in Section \ref{S:alg} and
we discuss the results obtained from it in Section \ref{S:results}.
Finally, we present an example of how the algorithm works for $\g$
of type $G_2$ and $e$ a short root vector in Section \ref{S:example}.

\section{Preliminaries}
\label{S:prelim}

We begin by giving the notation that we require, and recalling the definition
of the finite $W$-algebra $U(\g,e)$.
There are at present three equivalent definitions
of finite $W$-algebras in the literature. Here we only consider the
Whittaker model definition introduced in \cite{Pr2}:
this was proved to be equivalent to the definition
via BRST cohomology in \cite{DDDHK} and to the definition via Fedosov
quantization in \cite[Cor.\ 3.3.3]{Lo1}.

%Also we overlook the possible choice
%of a good grading to define the finite $W$-algebra: it is shown in
%\cite[Thm.\ 1]{BG} that the finite $W$-algebra is independent  of this choice
%up to isomorphism.  Finally, we note that we only consider the case
%where we choose a Lagrangian subspace $\l = \g(-1)^0$ of $\g(-1)$.
%An alternative equivalent definition valid for
%any isotropic subspace $\l \sub \g(-1)$ was
%given in \cite{GG}.

\subsection{}
Let $G$ be a simple simply connected algebraic group over $\bC$ and
let $\g$ be the Lie algebra of $\g$.  Let $e \in \g$ be a nilpotent
element and let $(e,h,f)$ be an $\sl_2$-triple in $\g$. Let
$(\cdot\,,\cdot)$ be a non-degenerate symmetric invariant bilinear
form on $\g$ normalized so that $(e,f)=1$.  Define $\chi \in \g^*$
by $\chi(x) = (e,x)$ for $x \in \g$.  The {\em Dynkin grading} $\g =
\bigoplus_{j \in \bZ} \g(j)$ of $\g$ is defined by $\g(j) = \{x \in
\g \mid [h,x] = jx\}$. We write $\g^e$ for the centralizer of $e$ in
$\g$.

Let $\t^e$ be a Cartan subalgebra of $\g^e \cap \g(0)$, and
let $\t$ be a Cartan subalgebra of $\g$ containing $\t^e$.  We write
$\Phi \sub \t^*$ for the root system of $\g$ with respect to $\t$
and $\Pi$ for a set of simple of roots. Recall that the {\em
restricted root system} $\Phi^e$ is defined by $\Phi^e = \{ \alpha
|_{\t^e} \mid \alpha \in \Phi\} \sub (\t^e)^*$, see \cite[\S 2 and
\S 3]{BG}. We have the $\t^e$-weight space decomposition $\g = \g_0
\oplus \bigoplus_{\alpha \in \Phi^e} \g_\alpha$ of $\g$, where
$\g_\alpha = \{x \in \g \mid [t,x] = \alpha(t)x \text{ for all } t
\in \t^e\}$ for $\alpha \in \Phi^e \cup \{0\}$.
%, so $\g_0$ is the centralizer of $\t^e$ in $\g$.
%As explained in
%\cite[\S 2]{BG}, we may choose a set of positive roots $\Phi_+e$ in
%$\Phi^e$; this is equivalent to choosing a parabolic subalgebra of
%$\g$ with Levi factor $\g_0$.

\subsection{}
For $x,y \in \g(-1)$, let $\langle x , y \rangle = \chi([x,y])$, so
that $\langle \cdot \,, \cdot \rangle$ defines a non-degenerate
alternating bilinear
form on $\g(-1)$.  Choose a Lagrangian subspace  $\g(-1)^0$ of $\g(-1)$
with respect to $\langle \cdot \,, \cdot \rangle $ and define the
nilpotent subalgebra
\[
\m = \g(-1)^0 \oplus \bigoplus_{i\leq -2}\g(i)
\]
of $\g$.
It is straightforward to check that $\chi$ restricts to a character of
$\m$, so we can consider the $1$-dimensional $U(\m)$-module $\bC_\chi$ and
the induced module $Q_\chi = U(\g) \otimes_{U(\m)} \bC_\chi$.  The {\em finite
$W$-algebra} associated to $\g$ and $e$
is defined to be the endomorphism algebra
\[
U(\g,e) = \End_{U(\g)}(Q_\chi)^\op.
\]

Let $I_\chi$ be the left ideal of $U(\g)$ generated by
all $x-\chi(x)$ for $x \in \m$.  The PBW-theorem for $U(\g)$
implies that $Q_\chi$ identifies with
$U(\g)/I_\chi$ as a vector space. Through this identification
we get an isomorphism
between $U(\g,e)$ and the space of $\m$-invariants of $U(\g)/I_\chi$:
\[
\{u + I_\chi \in U(\g) / I_\chi \mid [x,u] \in I_\chi \text{ for all }
x \in \m\}.
\]
We use these identifications throughout the sequel.

Thanks to \cite[Lem.\ 2.4]{Pr3}, there is an embedding $\t^e \into U(\g,e)$.
We identify $\t^e$ with its image in $U(\g,e)$, so we have an adjoint action
of $\t^e$ on $U(\g,e)$.
The weights of this action lie in $\bZ\Phi^e \sub (\t^e)^*$.

\subsection{} \label{ss:PBW}
In Theorem \ref{T:PBW} below we state a version of the PBW-theorem
for $U(\g,e)$.  To state this explicitly we need to pick a suitable
basis of $\g$.

Define the parabolic subalgebra $\p = \bigoplus_{j \ge 0} \g(j)$ of
$\g$; note that we have the inclusion $\g^e \sub \p$.  Let $x_1,
\dots,x_r$ be a basis of $\g^e$ and extend this to a basis
$x_1,\dots,x_m$ of $\p$; we assume that $x_i = e$ for some $i \in
\{1,\dots, r\}$.  The form $\langle \cdot \,, \cdot \rangle$ is
$\t^e$-invariant, i.e.\ $\langle [x,t],y \rangle = \langle x, [t,y]
\rangle$ for $x,y \in \g(-1)$ and $t \in \t^e$, so we may choose a
Witt basis $z_1,\dots,z_s,z^*_1,\dots,z^*_s$ of $\g(-1)$ consisting
of $\t^e$-weight vectors.  We choose $\g(-1)^0$ to have basis
$z^*_1,\dots,z^*_s$, and set $x_{m+i} = z_i$ and $x_{m+s+i} = z^*_i$
for $i = 1,\dots,s$. Next we choose a basis
$x_{m+2s+1},\dots,x_{m+2s+s'}$ of $\g(-2)$ as follows. It is clear
that $\ker \chi|_{\g(-2)}$ is $\t^e$-stable, so we choose
$x_{m+2s+1},\dots,x_{m+2s+s'-1}$ to be $\t^e$-weight vectors forming
a basis of $\ker \chi|_{\g(-2)}$ and $x_{m+2s+s'} = f$. Finally, we
extend to a basis $x_1,\dots,x_n$ of all of $\g$.

It is clear that we can choose all of the elements of our basis to
be weight vectors for $\t^e$ and eigenvectors for $\ad h$.  We let
$n_i \in \bZ$ and $\beta_i \in \Phi^e$ be such that  $x_i \in
\g(n_i) \cap \g_{\beta_i}$ for $i = 1,\dots,n$.

%with $x_i \in \g_{\beta_i}$ for $\beta_i \in \Phi^e \cup \{0\}$
%such that $x_i = f$ for some $i$ and $\chi(x_i) = 0$ for all other
%$i = m+2s+1,\dots,m+2s+l$.
%Finally, we extend this to a basis $x_1,\dots,x_n$ of all of $\g$ taking
%$x_i \in \g(n_i) \cap \g_{\beta_i}$ for $i = m+2s+l+1,\dots,n$, where
%$n_i \le -3$ and $\beta_i \in \Phi^e \cup \{0\}$.

%Since $\g^e$ and $\p$ are $\t^e$-stable, we may choose the $x_i$
%so that $x_i \in \g(n_i) \cap \g_{\beta_i}$ for some $n_i \in \bZ_{\ge 0}$
%and $\beta_i \in \Phi^e \cup \{0\}$.

%; we write $\beta_i \in \Phi^e$
%for the $\t^e$-weight of each $x_i$ noting that $\g_0 \cap \g(-1) = \{0\}$.

Let $\{e_\alpha \mid \alpha \in \Phi\} \cup \{h_\alpha \mid \alpha \in
\Pi\}$ be a Chevalley basis of $\g$.  Let $\g_\bZ$ be the corresponding
$\bZ$-form on $\g$ and let $\g_\bQ = \g_\bZ \otimes_\bZ \bQ$.  It is well
known that the $G$-orbit of $e$ intersects $\g_\bZ$, so we can assume that
$e \in \g_\bZ$.  Then it is easy to see that we can choose our basis
elements $x_1,\dots,x_n \in \g_\bQ$.

From the PBW-theorem for $U(\g)$ we see that a basis of $Q_\chi$ is
given by the cosets $x^\va + I_ \chi = x_1^{a_1}\cdots
x_{m+s}^{a_{m+s}} + I_\chi$ for $\va =(a_1,\dots,a_{m+s}) \in
\bZ_{\ge 0}^{m+s}$. For $\va \in \bZ_{\ge 0}^{m+s}$, we define
\[
|\va| = \sum_{i=1}^{m+s} a_i \quad \text{and} \quad
|\va|_e = \sum_{i=1}^{m+s} a_i(n_i + 2).
\]
This allows us to define the {\em Kazhdan filtration} on $Q_\chi$ by
declaring that $x^\va + I_\chi$ has filtered degree $|\va|_e$.
%; we write $\deg_e x^\va = |\va|_e$.
This restricts to the {\em Kazhdan filtration} on $U(\g,e)$ and we
write $F_i U(\g,e)$ for the $i$th filtered part of the Kazhdan
filtration on $U(\g,e)$.  Given an element $u \in U(\g,e)$, we say
that $u$ has Kazhdan degree $i$ to mean $u \in F_i U(\g,e)$, not
necessarily assuming that $i$ is minimal.

The PBW-theorem for $U(\g,e)$ is sometimes stated in terms of the
associated graded algebra $\gr U(\g,e)$ for the Kazhdan filtration,
see for example \cite[Thm.\ 4.1]{GG}.  For our purposes the
following more explicit version is more convenient.  The statement
combines \cite[Thm.~4.1]{Pr1} and \cite[Lem.\ 2.2]{Pr3}.  The
uniqueness statements are not given in the references, but they are
straightforward to deduce.

%\newpage

\begin{thm}
\label{T:PBW}
Let $x_1,\dots,x_n$ be a basis of $\g$ as above.
Then the following hold:
\begin{enumerate}
\item There is a set of generators for $U(\g,e)$ given by
\[
\Theta_i = \left( x_i + \sum_{|\va|_e
\leq n_i+2} \lambda^i_\va x^\va \right)+ I_\chi,
\]
for $i=1,\dots,r$, where the coefficients
$\lambda^i_\va \in \bQ$ are zero when
$a_{r+1} = \cdots = a_{m+s} = 0$, or if $|\va|_e = n_i + 2$
and $\vert \va \vert = 1$. The coefficients $\lambda^i_\va$
are uniquely determined by the choice of ordered basis
$x_1,\dots,x_n$ of $\g$ and the above vanishing conditions.
\item The $\Theta_i$ are weight vectors for $\t^e$ with weight
$\beta_i$.
\item
The monomials $\Theta^\va = \Theta_1^{a_1} \cdots \Theta_r^{a_r}$
with $\va \in \bZ_{\ge 0}^r$ form a PBW-basis of $U(\g,e)$.
\item
We have $[\Theta_i,\Theta_j] \in F_{n_i+n_j +2} U(\g,e)$ for
$i,j = 1,\dots,r$.  Moreover, if $[x_i,x_j] = \sum_{k=1}^r \mu_{ij}^k x_k$
in $\g^e$, then
\[
[\Theta_i,\Theta_j] = \sum_{k=1}^r \mu_{ij}^k \Theta_k +
q_{ij}(\Theta_1,\dots,\Theta_r) \mod F_{n_i+n_j} U(\g,e),
\]
where $q_{ij}$ is a polynomial with coefficients in $\bQ$, and zero
constant and linear terms.
\end{enumerate}
\end{thm}

%\subsection{} We finish the preliminaries by recalling some results of
%highest weight theory from \cite{BGK}.

%{\em Not sure what will go here yet.  Will do this quite late if we decide to
%give highest weights in table of results.}

\section{Removing relations} \label{S:small}

Theorem \ref{T:PBW} allows one to determine a presentation of
$U(\g,e)$, as in \cite[Lem.\ 4.1]{Pr4}, from which one can work out
all $1$-dimensional representations of $U(\g,e)$. However, this
presentation involves a large number of commutator relations, so is
rather laborious to calculate.  In Theorem \ref{T:small} we show
that fewer relations suffice to obtain a presentation of $U(\g,e)$.
Further in Proposition \ref{P:1-dim}, we show that if we are only
interested in determining the $1$-dimensional representations of
$U(\g,e)$, then we need to consider even fewer relations.

In the following lemma, which is required for the proof
of Theorem \ref{T:small}, we use the
notation introduced in \S \ref{ss:PBW} and Theorem \ref{T:PBW}.

\begin{lem} \label{L:gens}
Suppose that $\g^e$ is generated by $x_1,\dots,x_b$ for $b \le r$.  Then
$U(\g,e)$ is generated by $\Theta_1,\dots,\Theta_b$.
\end{lem}

\begin{proof}
We may assume that $x_{b+1},\dots,x_r$ are chosen so that $n_{b+1}
\le \dots \le n_r$.  The assumption that $x_1,\dots,x_b$ generate
$\g^e$ implies that there exists $i,j \le b$ such that $[x_i,x_j]$
is in $\g(n_{b+1})$, but does not lie in the subspace of $\g^e$
spanned by $x_1,\dots,x_b$.  Now we can assume that we picked our
basis of $\g^e$ with  $x_{b+1} = [x_i,x_j]$.  Similar arguments show
that for $k > b$ we can assume that
\[
x_k = \sum_{i,j < k} \nu_{ij}^k [x_i,x_j],
\]
where $\nu_{ij}^k \in \bQ$.

To prove the lemma, it suffices to show that $\Theta_i$ lies in the
subalgebra $W$ of $U(\g,e)$ generated by $\Theta_1,\dots,\Theta_b$
for $i = b+1,\dots,r$. Suppose that we have shown that
$\Theta_{b+1},\dots,\Theta_{k-1} \in W$; in particular, this means
that $F_{n_k+1} U(\g,e) \sub W$.   From Theorem \ref{T:PBW}(4) we
see that
\[
\sum_{i,j < k} \nu_{ij}^k [\Theta_i,\Theta_j] = \Theta_k +
G_k(\Theta_1,\dots,\Theta_{k-1}) + H_k(\Theta_1,\dots,\Theta_{k-1}),
\]
where $G_k$ and $H_k$ are a polynomials over $\bQ$ such that
$G_k(\Theta_1,\dots,\Theta_{k-1}) \in F_{n_k+2} U(\g,e)$,
$H_k(\Theta_1,\dots,\Theta_{k-1}) \in F_{n_k} U(\g,e)$, and $G_k$
has zero constant and linear terms. Therefore,
$G_k(\Theta_1,\dots,\Theta_{k-1})$ can be written as a sum of
products of elements of $ F_{n_k+1} U(\g,e)$ and
$H_k(\Theta_1,\dots,\Theta_{k-1}) \in F_{n_k+1} U(\g,e)$. Also each
commutator $[\Theta_i,\Theta_j] \in W$, so we have $\Theta_k \in W$,
as required.
\end{proof}

Although Lemma \ref{L:gens} shows that we can get by with fewer
generators, we in fact use it to show that some relations are not
needed for a presentation of $U(\g,e)$ in Theorem \ref{T:small}.
The case where $b = r$ in Theorem \ref{T:small} is \cite[Lem.\
4.1]{Pr4}. Throughout the proof and statement of Theorem
\ref{T:small} we use $(\vTheta)$ as a shorthand for
$(\Theta_1,\dots,\Theta_r)$.

\begin{thm} \label{T:small}
Suppose that $\g^e$ is generated by $x_1,\dots,x_b$ for $b \le r$.
Then $U(\g,e)$ is generated by $\Theta_1,\dots,\Theta_r$ subject only to
the relations
\[
[\Theta_i,\Theta_j] = F_{ij}(\Theta_1,\dots,\Theta_r) = F_{ij}(\vTheta),
\]
for $i = 1,\dots,b$ and $j=1,\dots,r$,
where $F_{ij}$ is a polynomial with coefficients in $\bQ$, and
$F_{ij}(\vTheta) \in F_{n_i+n_j+2} U(\g,e)$.
\end{thm}

\begin{proof}
By Theorem \ref{T:PBW}, the commutator $[\Theta_i,\Theta_j]$ is of the form
$F_{ij}(\vTheta)$, where $F_{ij}$ is a polynomial
satisfying the stated conditions, for $i,j = 1,\dots,r$.  Thanks to
\cite[Lem.\ 4.1]{Pr4}, $U(\g,e)$ is generated by the $\Theta_i$
subject only to these commutator relations.  Therefore, to
prove the theorem it suffices to show that the polynomials
$F_{kl}(\vTheta)$ for $k,l=b+1,\dots,r$
can be determined from the polynomials $F_{ij}(\vTheta)$
for $i = 1,\dots,b$ and $j=1,\dots,r$.

As in the proof of Lemma \ref{L:gens}, we assume that
$x_{b+1},\dots,x_r$ are chosen so that $n_{b+1} \le \dots \le n_r$
and that
\[
\Theta_k = \sum_{i,j < k} \nu_{ij}^k [\Theta_i,\Theta_j] -
G_k(\vTheta) - H_k(\vTheta),
\]
for $k > b$.
%, where we use the notation of Lemma \ref{L:gens}.

Consider $[\Theta_k,\Theta_l]$
% = F_{kl}(\Theta_1,\dots,\Theta_r)$
for some $k,l \in \{ b+1,\dots,r\}$.  We assume
inductively that we have calculated all $F_{k'l'}$
for $n_{k'} + n_{l'} < n_k + n_l$ in terms of the $F_{ij}$ for
$i=1,\dots,b$ and $j = 1,\dots,r$.
We have
\[
[\Theta_k,\Theta_l] = \sum_{i,j < k} \nu_{ij}^k [[\Theta_i,\Theta_j],\Theta_l]
- [G_k(\vTheta),\Theta_l] - [H_k(\vTheta),\Theta_l].
\]

Consider a term of the form $[[\Theta_i,\Theta_j],\Theta_l]$. This commutator
has Kazhdan degree $n_i+n_j+n_l+6$ (note that $n_k =n_i+n_j$).
We can apply the Jacobi identity to obtain
\begin{align*}
[[\Theta_i,\Theta_j],\Theta_l] &= [[\Theta_i,\Theta_l],\Theta_j] +
[\Theta_i,[\Theta_j,\Theta_l]] \\
&= [F_{il}(\vTheta),\Theta_j] + [\Theta_i,F_{jl}(\vTheta)].
\end{align*}
By induction, we can determine $F_{il}(\vTheta)$ and
$F_{jl}(\vTheta)$ and they both have Kazhdan degree $n_k+n_l+4$. Now
we can apply the Leibniz rule and inductive hypothesis to calculate
$[F_{il}(\vTheta),\Theta_j]$ and $[\Theta_i,F_{jl}(\vTheta)]$, as
polynomials in $\vTheta$ of Kazhdan degree $n_k+n_l+2$.

Next consider the term $[G_k(\vTheta),\Theta_l]$, this has Kazhdan
degree $n_k+n_l+4$, and $G_k(\vTheta)$ has zero linear term.
Therefore, we can apply the Leibniz rule and the inductive
hypothesis to calculate $[G_k(\vTheta),\Theta_l]$ as a polynomial in
$\vTheta$ of Kazhdan degree $n_k+n_l+2$.

Finally, consider the term $[H_k(\vTheta),\Theta_l]$.  Since
$H_k(\vTheta)$ has Kazhdan degree $n_k$, this can be calculated as a
polynomial in $\vTheta$ of Kazhdan degree $n_k+n_l$ using the
Leibniz rule and inductive hypothesis.

Thus we have determined $[\Theta_k,\Theta_l] = F_{kl}(\vTheta)$ as a polynomial
of Kazhdan degree $n_k+n_l+2$, as required.
\end{proof}

%\begin{rem}
%It is natural to ask whether it is possible to obtain presentations of
%$U(\g,e)$ with fewer generators and relations than those given in Theorem
%\ref{T:small}.  Of course, it is possible to use Lemma \ref{L:gens} to
%reduce the number of generators... probably make relations horrible.
%References to suggest that we can: \cite[Thm.\ 10.1]{BK1},
%\cite[Thm.\ 1.1]{Pr3}?.
%\end{rem}

We next discuss how to find all $1$-dimensional representations of
$U(\g,e)$ from a presentation as in Theorem \ref{T:small}. From now
we fix $b$ such that $x_1,\dots,x_b$ generate $\g^e$. Let $\rho :
U(\g,e) \to \bC$ be a $1$-dimensional representation of $U(\g,e)$.
Then $\rho$ is determined by the values $\rho(\Theta_i)$ for $i =
1,\dots,r$. These must satisfy the relations
$[\rho(\Theta_i),\rho(\Theta_j)] = F_{ij}(\rho(\Theta_1), \dots,
\rho(\Theta_r))$ for $i = 1,\dots,b$ and $j = 1,\dots,r$. We thus
see that finding 1-dimensional representations of $U(\g,e)$ is
equivalent to finding solutions to the polynomial equations
\begin{equation} \label{e:1dsols}
F_{ij}(t_1,\dots,t_r) = 0 \text{ for $i = 1,\dots,b$ and $j = 1,\dots,r$,}
\end{equation}
with $t_1,\dots,t_r \in \bC$.

The following lemma and proposition show that the values of certain $t_i$
in solutions to \eqref{e:1dsols} are forced to be zero, which means
that we do not have to consider all the equations in \eqref{e:1dsols} to
work out all $1$-dimensional representations of $U(\g,e)$. This in turn means
that if we only wish to determine the $1$-dimensional representations
of $U(\g,e)$, then it is not necessary to calculate all the commutator
relations from Theorem \ref{T:small}.

\begin{lem} \label{L:weights}
Let $\rho : U(\g,e) \to \bC$ be a representation of $U(\g,e)$.  Then
$\rho(\Theta_i) = 0 $ for all $i = 1,\dots,r$ such that
$\beta_i \ne 0$.
\end{lem}

\begin{proof}
As explained in \S 2.2, there is an embedding $\t^e \into U(\g,e)$ and we identify
$\t^e$ with its image in $U(\g,e)$.  Let $i \in \{1,\dots,r\}$ such that
$\beta_i \ne 0$ and let $t \in \t^e$ with $\beta_i(t) \ne 0$.  Then
\[
\beta_i(t)\rho(\Theta_i) = \rho(\beta_i(t)\Theta_i)
= \rho([t,\Theta_i]) = [\rho(t),\rho(\Theta_i)] =0.
\]
Thus $\rho(\Theta_i) = 0$, as required.
\end{proof}

The next proposition
means that in order to determine the $1$-dimensional representations of
$U(\g,e)$, we only need to know the commutators $[\Theta_i,\Theta_j]$, when
$\beta_j = -\beta_i$.  We need to introduce some notation in order to state
the proposition.

Let $I = \{i \in \{1,\dots,r\} \mid \beta_i = 0\}$ and let
$J = \{(i,j) \in \{1,\dots,b\} \times \{1,\dots,r\}
\mid \beta_j = -\beta_i\}$.
Let $F_{ij}$ be the polynomials from Theorem \ref{T:small}
and for $(j,k) \in J$ define the polynomials $\bar F_{jk}
\in \bQ[T_i \mid i \in I]$ by $\bar F_{jk}(T_i \mid i \in I) =
F_{jk}(\delta_1T_1,\dots,\delta_r T_r)$, where $\delta_i = 1$ if $i \in I$
and $\delta_i = 0$ if $i \notin I$.

\begin{prop} \label{P:1-dim}
The 1-dimensional representations $\rho : U(\g,e) \to \bC$ are in
bijective correspondence with solutions $(t_i \mid i \in I)$ to the
polynomial equations $\bar F_{jk}(t_i \mid i \in I) = 0$ for $(j,k)
\in J$.  The solution $(t_i \mid i \in I)$ corresponds to the
$1$-dimensional representation $\rho$ determined by $\rho(\Theta_i)
= t_i$ for $i \in I$ and $\rho(\Theta_i) = 0$ if $i \notin I$.
\end{prop}

\begin{proof}
As discussed above, the $1$-dimensional representations of $U(\g,e)$
are given by solutions to \eqref{e:1dsols}. By Lemma
\ref{L:weights}, we must have $t_i = 0$ for $i \notin I$ in any
solution to these equations.  By considering the $\t^e$-weights, we
see that each monomial in a polynomial $F_{jk}$ for $(j,k) \notin J$
must contain a $T_i$ for some $i \notin I$.   The result follows
form these two observations.
\end{proof}

\section{The algorithm}
\label{S:alg}

In this section we describe our algorithm for calculating a presentation
of $U(\g,e)$ as given in Theorem \ref{T:small}.  We work with
a Chevalley $\bZ$-form of $\g$, and end up with a presentation
of $U(\g,e)$ that is defined over $\bZ[d^{-1}]$, where $d$ is a product
of certain primes.  We wish to keep $d$ as small as possible, but at some
points it is necessary to allow divisions by certain primes, as in
the construction in \cite[\S 2]{Pr5}.  Keeping $d$ small allows
us to control the bound on $p$ in Theorem \ref{T:mindim},
see Remark \ref{R:char}.
We keep track of how $d$ ``grows'' throughout the algorithm
and begin with $d$ being the product of bad primes for $\g$.

\subsection{Input}
The input to our algorithm is the root system
$\Phi$ of a simple Lie algebra
$\g$ over $\bC$ along with a set of simple roots $\Pi$ and the
labelled Dynkin diagram $D$ of a nilpotent orbit
in $\g$, we write $d_\alpha$ for the label of $\alpha \in \Pi$
in $D$.

From this data we can determine the Chevalley $\bZ$-form $\g_\bZ$ of
$\g$ with Chevalley basis $\{h_\alpha \mid \alpha \in \Pi\} \cup
\{e_\alpha \mid \alpha \in \Phi\}$. The Cartan subalgebra $\t$ has
basis $\{h_\alpha \mid \alpha \in \Pi\}$.
%Let $\g_\bQ = \g_\bZ \otimes_\bZ \bQ$ and $\t_\bQ = \t_\bZ \otimes_\bZ \bQ$;
%further,
We define $\g_{\bZ[d^{-1}]} = \g_\bZ \otimes_\bZ \bZ[d^{-1}]$. The
labelled Dynkin diagram $D$ determines a decomposition $\Phi =
\bigcup_{j \in \bZ} \Phi(j)$, where $\Phi(j) = \{\sum_{\alpha \in
\Pi} a_\alpha \alpha \in \Phi \mid \sum_{\alpha \in \Pi} a_\alpha
d_\alpha = j\}$; this in turn gives the grading $\g = \bigoplus_{j
\in \bZ} \g(j)$, where $\g(j)$ is spanned by the $e_\alpha$ with
$\alpha \in \Phi(j)$ for $j \ne 0$, and $\g(0)$ is spanned by $\t$
and the $e_\alpha$ with $\alpha \in \Phi(0)$.

\subsection{Finding the $\sl_2$-triple} \label{ss:triple}
We wish to determine the $\sl_2$-triple $(e,h,f)$.  The labelled
Dynkin diagram $D$ uniquely determines $h \in \t$ as follows:
we may write $h = \sum_{\alpha \in \Pi} \lambda_\alpha
h_\alpha$, and the $\lambda_\alpha$ are uniquely determined by the conditions
$\sum_{\alpha \in \Pi} \lambda_\alpha \beta(h_\alpha) = d_\beta$ for all
$\beta \in \Pi$.  In fact, since the eigenvalues of $h$ on each fundamental
irreducible representation of $\g$ are integers, we
have $\lambda_\alpha \in \bZ_{\ge 0}$ for all $\alpha$.
We move on to determining $e$. It is well-known that we may
choose $e$ to be of the form
$e = \sum_{\alpha \in \Gamma} e_\alpha \in \g(2)$,
where $\Gamma$ is a certain linearly independent subset of $\Phi(2)$.
%: for example, this can be seen for the
%$\g$ of classical type from the descriptions in \cite[\S 4--6]{EK} and
%for $\g$ of exceptional types from the tables in \cite[App.\ A]{dG},
%it is also possible to give a direct proof of this without
%appealing to a case by case analysis.
Now there is a unique $f \in \g$ such that $(e,h,f)$ is an $\sl_2$-triple.
We can find $f$ by writing $f = \sum_{\alpha \in \Phi(2)} \mu_\alpha
e_{-\alpha}$, then solving for the $\mu_\alpha$ in $[e,f] = h$.  One can show
that it is always possible to choose $e$, so that we get $f \in \g_{\bZ}$:
to do this one reduces inductively to the case where $e$ does not
lie in any proper subalgebra of $\g$ containing $\t$, and makes an
explicit calculation in the remaining cases.

\subsection{Determining the bilinear form}
Now that we have our $\sl_2$-triple we can determine our bilinear form
$(\cdot \,, \cdot)$ on $\g$.
We obtain this by rescaling the Killing form $\kappa$
so that $(e,f) = 1$.  In general it will not be the case that $(x,y) \in \bZ$
for $x, y \in \g_\bZ$; we require that $(\cdot \,, \cdot) :
\g_{\bZ[d^{-1}]} \times \g_{\bZ[d^{-1}]} \to \bZ[d^{-1}]$.  Thus
we may increase $d$ here to ensure that $\kappa(e,f)$ is
invertible in $\bZ[d^{-1}]$.
With this choice of bilinear form we have that $\chi \in \g^*$ defined by
$\chi(x) = (x,e)$ is such that $\chi : \g_{\bZ[d^{-1}]} \to \bZ[d^{-1}]$.

\subsection{Determining the basis} \label{ss:basis}
We next wish to find a basis of $\g$ as described in \S \ref{ss:PBW}.

Determining a basis $x_1,\dots,x_r$ of $\g^e$ is straightforward,
using the Chevalley commutator relations, though we may need to
increase $d$ here so that it is a $\bZ[d^{-1}]$-basis of $\g^e \cap
\g_{\bZ[d^{-1}]}$. We order this basis so that $\g^e$ is generated
by $x_1,\dots,x_b$ and this generating set is chosen so that it is
small and convenient to work with. In particular, we want
$x_1,\dots,x_b$ to generate $\g^e \cap \g_{\bZ[d^{-1}]}$ over
$\bZ[d^{-1}]$, so there is a possibility that we may have to
increase $d$ further.

We extend this to get a basis $x_1,\dots,x_m$ of $\p$.
We may need to increase $d$ so that $x_1,\dots,x_m$ form a
$\bZ[d^{-1}]$-basis of $\p_{\bZ[d^{-1}]} = \g_{\bZ[d^{-1}]} \cap \p$.

Next we move on to determine the basis $x_{m+1},\dots,x_{m+2s}$ of
$\g(-1)$.  We can choose $x_{m+1},\dots,x_{m+s}$ by picking a set of
positive roots $\Phi_+^e$ in $\Phi^e$.  Then we can readily take
$\{x_{m+1},\dots,x_{m+s}\} = \{e_\alpha \mid \alpha \in \Phi(-1)
\text{ and } \alpha |_{\t^e} \in \Phi^e_+\}$.  We take each
$x_{m+s+i}$ to be a linear combination (over $\bQ$) of $\{e_\alpha
\mid \alpha \in \Phi(-1) \text{ and } \alpha |_{\t^e} \in
-\Phi^e_+\}$, so that they are uniquely determined subject to
$x_{m+1},\dots,x_{m+2s}$ being a Witt basis of $\g(-1)$ with respect
to $\langle \cdot \,,\cdot \rangle$.  We may have to increase $d$ so
that all the coefficients that occur in these linear combinations
lie in $\bZ[d^{-1}]$.

To determine our basis $x_{m+2s+1},\dots,x_{m+2s+s'}$
of $\g(-2)$, we first choose a basis of the kernel
of $\chi$ restricted to $\g(-2)$ to consist of the $e_{-\alpha}$ for
$\alpha \in \Phi(2) \setminus \Gamma$, and some elements of the form
$e_{-\alpha} + \lambda e_{-\beta}$ for $\alpha,\beta \in \Gamma$ and $\lambda
\in \bQ$.  We extend to a basis  of $\g(-2)$ by
adding $f$.  We may increase $d$ here
%so that it is divisible by
%any primes dividing coefficients $\lambda$
to ensure that $\lambda$ lies in
$\bZ[d^{-1}]$, thus to ensure that we obtain a
$\bZ[d^{-1}]$-basis of $\g_{\bZ[d^{-1}]} \cap \g(-2)$.

Finally, we extend to a basis $x_1,\dots,x_n$ of all of $\g$
by taking the $x_i$ for
$i = m{+}2s{+}s'{+}1,\dots,n$ to be the Chevalley basis elements
$e_\alpha$ for $\alpha \in \bigcup_{j \le -3} \Phi(j)$.

\begin{rem} \label{R:dsmaller}
Although there are many places above, where it seems necessary to increase
$d$, the cases that we have considered suggest that it suffices to
take $d$ to be the product of bad primes for $\g$ along with those dividing
$\kappa(e,f)$.  For example, it seems likely that it is always
possible to choose
$x_1,\dots,x_m$ so that $x_1,\dots,x_r$ is a $\bZ$-basis of
$\g_\bZ \cap \g^e$, and $x_1,\dots,x_m$ is a $\bZ$-basis of $\g_\bZ \cap \p$.
\end{rem}

\subsection{Finding generators} \label{ss:gens}
We wish to determine expressions for the $\Theta_i$ in the form
given in Theorem \ref{T:PBW}.   First we list the monomials $x^\va =
x_1^{a_1}\cdots x_{m+s}^{a_{m+s}}$ for which the coefficient
$\lambda_\va^i$ may be non-zero: according to Theorem \ref{T:PBW}
these are those satisfying $\sum_{j=1}^{m+s} a_j \beta_j = \beta_i$
and $\sum_{j=1}^{m+s} a_j(n_j+2) \le n_i+2$, and excluding those for
which either $a_{r+1} =\cdots =a_{m+s} =0$, (i.e.\ terms in
$U(\g^e)$), or $\sum_{j=1}^{m+s} a_j(n_j+2) = n_i+2$ and $|\va| = 1$
(i.e.\ single terms of equal Kazhdan degree to the leading term
$x_i$).

We choose $K \sub \{m+s+1,\dots,n\}$ so that the nilpotent
subalgebra $\m$ is generated by $\{x_j \mid j \in K\}$; we want to
pick $K$ to be small to reduce the amount of computation required.
For each $x_j$ with $j \in K$, we require $[x_j,\Theta_i] \in
I_\chi$ for $i = 1, \ldots, r$.
This condition can be calculated using the Chevalley
commutator relations and gives rise to a set of linear equations for
the coefficients $\lambda_\va^i$.  Solving these equations yields a
unique set of coefficients $\lambda_\va^i$ for $\Theta_i$.

In some cases, when $\Theta_i$ has large Kazhdan degree, the number
of linear equations that we are required to solve is too large for the
computation to be feasible. In these cases we try to construct
$\Theta_i$ in terms of the $\Theta_j$ for $j < i$. We chose our
basis of $\g$ so that $\g^e$ is generated by $x_1,\dots,x_b$, which
allows us to make this construction if $i > b$, where we follow the
idea of the proof of Lemma \ref{L:gens}.  In practice, we find an
element in $U(\g,e)$ in terms of the $\Theta_j$ with $j < i$ with
leading term $x_i$ and then subtract monomials in the $\Theta_j$ for
$j < i$, until the conditions of Theorem \ref{T:PBW} are met.
%In particular, the
%procedure for finding the relations requires that the only term
%in the generator which comes from the centralizer of $e$ is the
%leading term. Therefore, given some generating subset of (basis elements of)
%$\g^e$, it is necessary only to calculate the generator $\Theta_k$ by means
%of finding the coefficients $\lambda_\va^k$ where $x_k$ lies
%in that set.

The expressions for the generators $\Theta_i$ involve coefficients
$\lambda_\va^i \in \bQ$.  At this point we may need to increase $d$
to ensure that all these coefficients lie in $\bZ[d^{-1}]$.

\subsection{Finding relations} \label{ss:relns}
In order to determine the presentation of $U(\g,e)$ from Theorem
\ref{T:small}, we now just have to find the relations.
So we have to find polynomials $F_{ij}$ for $i = 1,\dots,b$ and
$j = 1,\dots,r$.  To find $F_{ij}$ we first evaluate
$[\Theta_i,\Theta_j]$ as an element of $U(\g)$, i.e.\ we
take our expressions $\Theta_i = u_i + I_\chi$ and
$\Theta_j = u_j + I_\chi$, with $u_i,u_j \in U(\g)$
and calculate $[u_i,u_j] + I_\chi$.
This determines an expression of the form
\[
[\Theta_i,\Theta_j] =
\sum_{|\va|_e \le n_i + n_j + 2} \mu_\va^{i,j} x^\va + I_\chi,
\]
with $\mu_\va^{i,j} \in \bQ$.  Amongst the $\va$ with $\mu_\va^{i,j} \ne 0$ and
$|\va|_e$ maximal, there must be one with $a_k = 0$ for all $k>r$.  Let
$\va_1,\dots,\va_c$ be all such $\va$.  Then we consider
\[
[\Theta_i,\Theta_j] - \sum_{k=1}^c \mu_{\va_k}^{i,j} \Theta^{\va_k}+
I_\chi.
\]
This is an element of $U(\g,e)$ and by construction it must have lower
Kazhdan degree than $[\Theta_i,\Theta_j]$.  We continue by subtracting
terms of maximal Kazhdan degree and after a finite number of steps we
obtain the required polynomial $F_{ij}$.

In order to have our presentation defined over $\bZ[d^{-1}]$, we may
have to increase $d$, so that all coefficients of the polynomials
$F_{ij}$ lie in $\bZ[d^{-1}]$.

\subsection{Determining 1-dimensional representations}
Once we have the presentation of $U(\g,e)$, we can determine all
1-dimensional representations of $U(\g,e)$ by solving the equations \eqref{e:1dsols}
or just those from Proposition \ref{P:1-dim}.  This is achieved using
standard techniques for solving polynomial equations.

\subsection{Implementation in GAP}
We have implemented this algorithm in the computer algebra language
GAP \cite{GAP} as explained below.

The Lie algebra $\g$ is created in GAP by taking the inbuilt
Chevalley basis of the simple Lie algebra over $\bQ$ and
constructing the required basis as explained in \S 4.4. By taking
the Lie products in GAP of these elements, we create a table of
structure constants, from which the function
\texttt{LieAlgebraByStructureConstants} returns $\g$ with the
required ordered basis, and the universal enveloping algebra $U(\g)$
is created using the function \texttt{UniversalEnvelopingAlgebra}.
This allows us to make all of the calculations required in the
algorithm.

It should be noted that for calculations and operations in the
universal enveloping algebra $U(\g)$ involving elements with many
terms  GAP functions can be particularly slow and require a lot of
memory. Such calculations can be speeded up by storing elements of
$U(\g)$ as elements of a particular polynomial ring rather than as
elements of the universal enveloping algebra.

% greatly reduces the memory
%required by GAP for these elements, and also speeds up any calculations
%with these terms. For multilpication of elements of $U(\g)$, the
%calculation is carried out term by term (as elements of $U(\g)$,
%then converted to elements of the polynomial ring), and taking
%the sum, making the calculation much faster.

\section{Results} \label{S:results}

We have used our algorithm to calculate presentations of $U(\g,e)$
for all cases where $\g$ is of type $G_2,F_4$ or $E_6$ and $e \in \g$
is rigid nilpotent.  From these
presentations we have determined all $1$-dimensional representations
of $U(\g,e)$.  For $\g$ of type $E_7$ and $e$ rigid nilpotent, we have
calculated all generators of $U(\g,e)$ and enough relations
to determine all $1$-dimensional representations
of $U(\g,e)$ using Proposition \ref{P:1-dim}.

In these cases we found that there are one or two $1$-dimensional
representations of $U(\g,e)$, as shown in Table \ref{Tab:rigidres}.
Here we give the Bala--Carter label of the rigid nilpotent orbits in
$\g$; the column says whether the number of $1$-dimensional
representations of $U(\g,e)$ is one or two.  A list of rigid orbits
can be found in \cite[p.\ 173]{Sp}.

\begin{table}[h!tb]
\renewcommand{\arraystretch}{1.5}
\begin{tabular}{|l|l|l|}
\hline
$\g$ & 1 & 2 \\
\hline
\hline
$G_2$ & $A_1$ & $\widetilde A_1$ \\
\hline
$F_4$ & $A_1$, \:\: $\widetilde A_1$, \:\: $A_1 + \widetilde A_1$, \:\:
$A_2 + \widetilde A_1$ &  $\widetilde A_2 + A_1$ \\
\hline
$E_6$ & $A_1$, \:\: $3A_1$, \:\: $2A_2 + A_1$ & \\
\hline
$E_7$ & $A_1$, \:\: $2A_1$, \:\: $(3A_1)'$, \:\: $4A_1$, \:\:
$A_2 + 2A_1$, \:\: $2A_2 + A_1$
& $(A_3 + A_1)'$ \\
\hline
\end{tabular}
\medskip
\caption{The number of $1$-dimensional representations of $U(\g,e)$
for $e$ rigid} \label{Tab:rigidres}
\end{table}

%In Table \ref{Tab:rigidres}, we record the
%results of our calculations.   Also
%include some $E_8$ cases, see Remark \ref{R:E8}.
%List of rigid orbits can be found in \cite[p 173]{Sp}.

%Give number $n$ of $1$-dimensionals
%Highest weight $\lambda$ is given;
%and value on $\Theta_i$s for $\beta_i$, others all go to zero.

%\begin{table}[h!tb]
%\renewcommand{\arraystretch}{1.5}
%\begin{tabular}{|l|l|l|l|l|}
%\hline $G$ & $G \cdot e$ & $n$ & $\Theta_i$s & $\lambda$ \\
%\hline\hline
%$G_2$ & $A_1$ & & & \\
%& $\widetilde A_1$ & 2 & & \\
%\hline\hline
%$F_4$ & $A_1$ & & & \\
%& $\widetilde A_1$ & & & \\
%& $A_1 + \widetilde A_1$ & & & \\
%& $A_2 + \widetilde A_1$ & & & \\
%& $\widetilde A_2 + A_1$ & 2 & & \\
%\hline
%\hline
%$E_6$ & $A_1$ & & & \\
%& $3 A_1$ & & & \\
%& $2A_2 + A_1$ & & & \\
%\hline
%\hline
%$E_7$ & $A_1$ & & & \\
%& $2A_1$ & & & \\
%& $(3A_1)'$ & & & \\
%& $4A_1$ & & & \\
%& $A_2 + 2A_1$ & & & \\
%& $2A_2 + A_1$ & & &\\
%& $(A_3 + A_1)'$ & 2 & (2,3,4,5,6,7,8) & (1,1,1,1,1,1,1) \\
%\hline
%\hline
%$E_8$ & ??? & & & \\
%\hline
%\end{tabular}
%\medskip
%\caption{The $1$-dimensional representations of $U(\g,e)$ for $e$ rigid} \label{Tab:rigidres}
%\end{table}

\begin{rem} \label{R:E8}
For some cases where $\g$ is of type $E_8$ and $e \in \g$ is rigid
nilpotent we are able to show that there exist $1$-dimensional
representations of $U(\g,e)$.  We have checked this when
$e$ has Bala--Carter label $A_1$, $2A_1$, $3A_1$, $4A_1$, $A_2 +
A_1$, $A_2 + 2A_1$, $A_2 + 3A_1$, $2A_2 + A_1$ and $A_3 + A_1$,
where in each but the last case there is just one $1$-dimensional
representation while in the last instance there are two such representations.
At present it is
computationally unfeasible to deal with the remaining 8
rigid nilpotent $e \in \g$.

Recall that the height of $e$ is defined to be the
maximal $j$ for which $\g(j) \ne 0$.  When the height of $e$ is large
then there are generators $\Theta_i$ of $U(\g,e)$ with large Kazhdan
degree.  Consequently, the expression for $\Theta_i$ given in
Theorem \ref{T:PBW} can be very complicated.  This means that it is
at present unfeasible to determine all generators, and to calculate
the required commutators.
\end{rem}

\begin{rem} \label{R:char}
We discuss some rationality issues related to Theorem
\ref{T:mindim}.  In our algorithm we keep track of a positive
integer $d$ such that $U(\g,e)$ is defined over $\bZ[d^{-1}]$.  For
$\g$ of type $G_2$, $F_4$ or $E_6$ and $e$ rigid, where we have
calculated a full presentation of $U(\g,e)$, we see that we can take
$d$ to be the product of the bad primes for $\g$ along with those
dividing $\kappa(e,f)$. This means that the reduction modulo $p$
argument from \cite[\S 2]{Pr5} goes through for all good $p$ not
dividing $\kappa(e,f)$.  This allows one to obtain an explicit lower
bound $M$ on $p$ for Theorem \ref{T:mindim} for $\g$ of type $G_2$,
$F_4$ or $E_6$, i.e.\ so that the conclusions of the theorem hold
for all $p>M$. This bound $M$ is determined by considering all pairs
$(\l,e_0)$, where $\l$ is a Levi subalgebra of $\g$ and $e_0 \in \l$
is rigid nilpotent, then taking $M$ to be the maximum prime dividing
$\kappa_\l(e_0,f_0)$, where $(e_0,h_0,f_0)$ is an $\sl_2$-triple in
$\l$ and $\kappa_\l$ is the Killing form on $\l$. To ensure that
this is the correct bound it was necessary to calculate
presentations of the finite $W$-algebras $U(\l,e_0)$ associated to
all such pairs $(\l,e_0)$. Explicitly, for $G_2$ we obtain $M=3$,
and for $F_4$ and $E_6$ we get $M=5$. We have to omit the prime $5$
for $\g$ of type $F_4$ or $E_6$ as it divides $\kappa(e_0,f_0)$,
where $e_0$ is a certain rigid nilpotent element for $\l$ of type
$B_3$ or $D_5$ respectively.

For the cases where $\g$ is of type $E_7$ it is not possible to be so explicit
about a bound.  We have not calculated all the relations, and we cannot rule
out the possibility that the other relations will lead to an increase in $d$.
In all the relations that we have calculated the only primes occurring
in denominators of coefficients are bad primes for $\g$, so it seems likely
that this is the case for all relations.  If this were true, then we would
get an analogous bound as for $\g$ in the previous paragraph.
\end{rem}

\newpage

\section{An example} \label{S:example}

We illustrate our algorithm with an explicit example
where $\g$ is of type $G_2$ and $e$ is a short root vector.

Let $\g$ be the simple Lie algebra of type $G_2$. Then GAP gives the
Chevalley basis $b_1, \dots ,b_{14}$, where $b_1, \dots ,b_6$ denote the
positive root vectors, $b_7, \dots, b_{12}$ denote the negative root
vectors, and the 2-dimensional Cartan subalgebra is generated by
$b_{13}=[b_1,b_7]$ and $b_{14}=[b_2,b_8]$.  We divide the Killing
form by $24$ to obtain our bilinear form $(\cdot \,, \cdot)$.

An $\sl_2$-triple for the orbit is determined to be $e=b_4$,
$h=2b_{13}+3b_{14}$ and $f=b_{10}$.  We construct our basis $x_1,
\dots ,x_{14}$ in terms of the Chevalley basis $b_1,\dots,b_{14}$ as
described in \S \ref{ss:basis}.  In Table \ref{Tab:g2} below, we
give this basis and show which parts of the basis form bases of
$\g^e$, $\p$ and $\m$, respectively.  We also give the values of $n_i$ and
$\beta_i$; as $\t^e = \bC x_{6}$ is $1$-dimensional, the weight
$\beta_i$ can be identified with the integer such that $[x_6,x_i] =
\beta_i x_i$.  We note that we can pick $x_1,\dots,x_m$ to be
multiples of elements of the Chevalley basis; this is not the case
in general.  We observe that our basis is a $\bZ[\frac{1}{2}]$-basis
%$\bZ[2^{-1}]$-basis
of
$\g_{\bZ[\frac{1}{2}]}$;
%$\g_{\bZ[2^{-1}]}$;
however, we only view it as a
$\bZ[\frac{1}{6}]$-basis
%$\bZ[6^{-1}]$-basis
of
$\g_{\bZ[\frac{1}{6}]}$,
%$\g_{\bZ[6^{-1}]}$,
as we have already divided by 3 to make $(e,f)=1$ (also $3$ is a bad
prime for $G_2$). We note that there is no small generating set of
$\g^e$, so we take $b=r=6$; also there is no small generating set of
$\m$, so we take $K$ from \S \ref{ss:gens} to be $\{11,12,13,14\}$.

\begin{table}[h!tb]
\begin{center}
\renewcommand{\arraystretch}{1.3}
\begin{tabular}[h]{|l|rrrrrrrrrrrrrr|}
\hline
& \multicolumn{9}{c}{$\p$} \vline & \multicolumn{1}{c}{} \vline &
\multicolumn{4}{c}{$\m$} \vline \\
\hline
& \multicolumn{6}{c}{$\g^e$} \vline & \multicolumn{8}{c}{}  \vline \\
\hline
 & $x_1$ & $x_2$ & $x_3$ & $x_4$ & $x_5$ & $x_6$ & $x_7$ & $x_8$
& $x_9$ & $x_{10}$ & $x_{11}$ & $x_{12}$ & $x_{13}$ & $x_{14}$ \\
\hline
 & $b_6$ & $b_5$ & $b_4$ & $b_2$ & $b_8$ & $b_{14}$ & $b_1$ & $b_3$
& $b_{13}$ & $b_9$ & $\frac{1}{2}b_7$ & $b_{10}$
& $b_{11}$ & $b_{12}$\\
\hline
$n_i$ & $3$ & $3$ & $2$ & $0$ & $0$ & $0$ & $1$ & $1$ & $0$
& $-1$ & $-1$ & $-2$ & $-3$ & $-3$ \\
\hline
$\beta_i$ & $1$ & $-1$ & $0$ & $2$ & $-2$ & 0 & $-1$ & $1$
& $0$ & $-1$ & $1$ & $0$ & $1$ & $-1$ \\
\hline
\end{tabular}
\end{center}
\medskip
\caption{Basis of $\g$ of type $G_2$}
\label{Tab:g2}
\end{table}

%\begin{center}
%\begin{tabular}[h]{|l|cccc|}
%\hline
%$x_i$ & $b_1$ & $b_3$ & $b_{13}$ & $b_9=z_1$ \\
%\hline
%$n_i$ & 1 & 1 & 0 & $-1$ \\
%\hline
%\end{tabular}
%\end{center}

%\begin{center}
%\begin{tabular}[h]{|l|cccc|}
%\hline
%& \multicolumn{4}{c}{$\m$} \vline\\
%\hline
%$x_i$ & $\frac{1}{2}b_7=z_1^\prime$ & $b_{10}=f$ & $b_{11}$ & $b_{12}$ \\
%\hline
%$n_i$ & $-1$ & $-2$ & $-3$ & $-3$ \\
%\hline
%$\chi(x_i)$ & $0$ & $1$ & $0$ & $0$ \\
%\hline
%\end{tabular}
%\end{center}
%\medskip}

We calculate $\Theta_i$ as explained in \S \ref{ss:gens}.  We
illustrate the procedure with the calculation of $\Theta_1$.  We
have that $\Theta_1$ is of the form $\Theta_1 = (x_1 + \sum
\lambda_\va^1 x^\va ) + I_\chi$ and satisfies the conditions of
Theorem \ref{T:PBW}. We first determine for which $\va$ the
coefficient $\lambda_\va^1$ can be non-zero and denote this set by
$A_1$; the elements of $A_1$ correspond to the monomials
$x_4x_6x_{10}$, $x_4x_7$, $x_4x_9x_{10}$, $x_4x_{10}$, $x_6x_8$,
$x_8$, and $x_8x_9$.  We next calculate $[x_i,x_1 + \sum_{\va \in
A_1} \lambda_\va^1 x^\va]$ in $U(\g)$ for $i = 11,\dots,14$, viewing
the $\lambda_\va^1$ as indeterminates.  Then we project  into
$Q_\chi = U(\g)/I_\chi$ to obtain a set of linear equations to solve
for the coefficients $\lambda_\va^1$.  These equations have a unique
solution which determines the value of $\Theta_1$.  We determine
$\Theta_i$ for $i = 2,\dots,6$ in the same way and obtain the
complete list of generators of $U(\g,e)$ given below:

\newpage

\begin{eqnarray*}
\Theta_1 &=& (x_1 + 3x_4x_6x_{10} + x_4x_7 + 2x_4x_9x_{10} -4x_4x_{10} \\
         & & +2x_6x_8 -4x_8 +x_8x_9) + I_\chi \\ & & \\
\Theta_2 &=& (x_2 + \tfrac{1}{2}x_3x_{10} -\tfrac{3}{2}x_4x_5x_{10} +
\tfrac{1}{2}x_4x_{10}^3 -x_5x_8 -x_6x_7\\
         & &  -\tfrac{1}{2}x_6x_9x_{10} + x_6x_{10} + 2x_7 -x_7x_9 +
\tfrac{1}{2}x_8x_{10}^2 \\
         & & + \tfrac{5}{2}x_9x_{10} -\tfrac{1}{2}x_9^2x_{10} -3x_{10}) +
I_\chi \\
&&\\
\Theta_3 &=& (x_3 + \tfrac{3}{4} x_4x_{10}^2 +3x_6x_9 +x_8x_{10} -5x_9 +
x_9^2) +I_\chi \\
&&\\
\Theta_4 &=& x_4 + I_\chi \\
&&\\
\Theta_5 &=& (x_5 - \tfrac{1}{4}x_{10}^2 ) + I_\chi \\
&&\\
\Theta_6 &=& x_6 + I_\chi.
\end{eqnarray*}

Next we calculate the relations. We illustrate this
by calculating the polynomial $F_{3,5}$ such that
$[\Theta_3,\Theta_5] = F_{3,5}(\Theta_1,\dots,\Theta_6)$.
We begin by taking the expressions for $\Theta_3$ and
$\Theta_5$ above and calculating
\[
[\Theta_3,\Theta_5] = (9x_5x_6 -\tfrac{9}{4}x_6x_{10}^2 -\tfrac{51}{2}x_5
+ \tfrac{15}{8}x_{10}^2) + I_\chi.
\]
The monomial $x_5x_6$ is the only monomial of Kazhdan degree
$n_3+n_5+2$ and consists of basis elements of $\g^e$ that occur
in the above expression.  Therefore, according to the method described in
\S \ref{ss:relns}, we consider
\[
[\Theta_3,\Theta_5] - 9\Theta_5 \Theta_6 = (-\tfrac{51}{2}x_5 -
\tfrac{51}{8} x_{10}^2 ) +I_\chi.
\]
Following the algorithm leads us to calculate
\[
[\Theta_3,\Theta_5] - 9\Theta_5 \Theta_6 + \tfrac{51}{2}\Theta_5 = 0 + I_\chi.
\]
Therefore, we have calculated the relation and we have that
$F_{3,5}(T_1,\dots,T_6) = 9T_5T_6 - \frac{51}{2}T_5$.

\newpage

The other commutator relations are found in the same way and all
are given below;
we omit all commutators that are equal to 0:
\begin{eqnarray*}
\lbrack \Theta_1,\Theta_2 \rbrack &=& 5\Theta_3\Theta_4\Theta_5-
\tfrac{1}{2}\Theta_3^2-\Theta_3\Theta_6^2+9\Theta_4\Theta_5\Theta_6^2 \\
&& -\tfrac{9}{2}\Theta_4^2\Theta_5^2+7\Theta_3\Theta_6-
\tfrac{69}{2}\Theta_4\Theta_5\Theta_6+6\Theta_6^3 \\
                    & & -6\Theta_3 + \tfrac{93}{4}\Theta_4\Theta_5-
30\Theta_6^2+42\Theta_6-18 \\
%&&\\
\lbrack \Theta_1,\Theta_3\rbrack  & = & 6\Theta_1\Theta_6-
3\Theta_2\Theta_4-3\Theta_1\\
\lbrack \Theta_1,\Theta_5\rbrack  & = & \Theta_2\\
\lbrack \Theta_1,\Theta_6\rbrack  & = & -\Theta_1\\
\lbrack \Theta_2,\Theta_3\rbrack  & = & -3\Theta_1\Theta_5-
6\Theta_2\Theta_6+12\Theta_2\\
\lbrack \Theta_2,\Theta_4\rbrack  & = & \Theta_1\\
\lbrack \Theta_2,\Theta_6\rbrack  & = & \Theta_2\\
\lbrack \Theta_3,\Theta_4\rbrack  & = & {}-9\Theta_4\Theta_6+
\tfrac{15}{2}\Theta_4\\
\lbrack \Theta_3,\Theta_5 \rbrack  & = &
9\Theta_5 \Theta_6 - \tfrac{51}{2}\Theta_5 \\
\lbrack \Theta_4,\Theta_5\rbrack  & = & \tfrac{1}{2} + \Theta_6\\
\lbrack \Theta_4,\Theta_6\rbrack  & = & -2\Theta_4\\
\lbrack \Theta_5,\Theta_6\rbrack  & = & 2\Theta_5.
\end{eqnarray*}
%with all the other commutators satisfying
%$[\Theta_i,\Theta_j] = 0$ for $1 \leq i,j \leq 6$.

We can observe that $2$ is the only prime
occurring in the denominators in the formulas for the generators and
relations.  Therefore, the finite $W$-algebra $U(\g,e)$ is defined over
$\bZ[\frac{1}{6}]$.
%\bZ[6^{-1}]$.

We move on to determining the 1-dimensional representations of $U(\g,e)$.
That is we have to find the solutions to the equations \eqref{e:1dsols}
where the polynomials $F_{ij}$ are given in the commutator
relations above.  We immediately see that we have
\[
t_1 = t_2 = t_4 = t_5 = 0 \text{ and } t_6 = -\tfrac{1}{2}.
\]
We substitute these values in $F_{1,2}$ and are left to solve
\[
-\tfrac{1}{2}t_3^2 - \tfrac{39}{4}t_3 - \tfrac{189}{4} = 0
\]
for $t_3$.  This gives the solutions $t_3 = -9,-\frac{21}{2}$, so there are two
$1$-dimensional representations of $U(\g,e)$.

\section*{Acknowledgments}

We are grateful to A.\ Premet for suggesting the problem considered
in this paper and for various helpful conversations. The last author
acknowledges the financial support of the EPSRC. Part of this paper
was written during a stay of the authors at the Isaac Newton
Institute for Mathematical Sciences, Cambridge during the
``Algebraic Lie Theory'' Programme in 2009.  We thank the referee
for some useful comments.

\end{document}